\documentclass[12pt]{article}
\usepackage{amsmath,amsfonts,amssymb,amsthm}

\numberwithin{equation}{section}

\newcommand{\Z}{{\mathbb Z}}

\newcommand{\C}{{\mathbb C}}

\newcommand{\al}{\alpha}

\renewcommand{\qed}{\mbox{ $\square$}}

\newtheorem{thm}{Theorem}[section]

\newtheorem{lem}[thm]{Lemma}

\theoremstyle{definition}
\newtheorem{remark}[thm]{Remark}

\begin{document}

\begin{center}
{\Large {\bfseries Decomposition of the vertex operator algebra 
$V_{\sqrt{2}D_l}$}}

\vspace{10mm}

Chongying Dong\footnote{Supported by NSF grant 
DMS-9700923 and a research grant from the Committee on Research, UC Santa Cruz.}\\[0pt]
Department of Mathematics, University of California\\[0pt]
Santa Cruz, CA 95064\\[0pt]
\vspace{3mm}
Ching Hung Lam\\[0pt]
Institute of Mathematics, University of
Tsukuba\\[0pt]
Tsukuba 305-8571, Japan\\[0pt]
\vspace{3mm}
Hiromichi Yamada\\[0pt]
Department of Mathematics, Hitotsubashi
University\\[0pt]
Kunitachi, Tokyo 186-8601, Japan
\end{center}

\vspace{10mm}

%\begin{abstract}
%\end{abstract}

\section{Introduction}

A weight two vector $v$ of a vertex operator algebra is called a
conformal vector with central charge $c$ if the component operators
$L_{v}(n)$ for $n\in {\mathbb Z}$ of $Y(v,z)=\sum_{n\in \mathbb
Z}L_v(n)z^{-n-2}$ satisfy the Virasoro algebra relation with central
charge $c$. In this case, the vertex operator subalgebra
$\mbox{Vir}(v)$ generated by $v$ is isomorphic to a
Virasoro vertex operator algebra with central charge $c$ (\cite{FZ},
\cite{M}).

Let $V_{\sqrt{2}R}$ be the vertex operator algebra associated with
$\sqrt{2}$ times an ordinary root lattice $R$ of type $A_{l}$,
$D_{l}$, or $E_{l}$. In \cite{DLMN}, several sets of mutually
orthogonal conformal vectors in $V_{\sqrt{2}R}$ were constructed
and studied. It was shown that the Virasoro element of
$V_{\sqrt{2}R}$ can be written as a sum of $l+1$ mutually
orthogonal conformal vectors $\omega ^{i}$, $1\le i\le l+1$. Since
these conformal vectors are mutually orthogonal, the subalgebra
$T$
generated by them is a tensor product $T=\otimes _{i=1}^{l+1}\mbox{Vir}
(\omega ^{i})$ and $\mbox{Vir}(\omega ^{i})$ is isomorphic to
$L(c_{i},0)$ with $c_{i}$ the central charge of $\omega ^{i}.$
Here, we denote by $L(c,h)$ the irreducible
highest weight module for the Virasoro
algebra with highest weight $h\in\C$ and central charge $c.$ 
Moreover, $V_{\sqrt{2}R}$ is completely reducible as a $T$-module
and each irreducible direct summand is in the form $\otimes
_{i=1}^{l+1}L(c_{i},h_{i})$. Hence one can study the structure of
the vertex operator algebra $V_{\sqrt{2}R}$ from a point of view
that $V_{\sqrt{2}R}$ is a $T$-module. Along this line, the
decomposition of $V_{\sqrt{2}A_{2}}$ and $V_{\sqrt{2}A_{3}}$ into
a direct sum of irreducible $T$-modules has been determined in
\cite{KMY} and \cite {DLY}. See also \cite{DMZ}, in which the
moonshine module $V^{\natural }$ is studied as a module for
$L(\frac{1}{2},0)^{\otimes 48}$.

In this paper we determine the decomposition of 
$V_{\sqrt{2}D_{l}}$ into a sum of irreducible $T$-modules for general
$l.$ There are many ways to choose a set of the mutually
orthogonal conformal vectors in \cite{DLMN} (see  Remark \ref{choice}). 
In this paper, we shall choose a set of conformal vectors so that
$c_{1}=1/2$, $c_{2}=7/10$, $c_{3}=4/5$, and $c_{i}=1$ for $4\le i\le l+1.$
Such a choice enables us to use the ideas and methods 
for the decompostion of $V_{\sqrt{2}A_{3}}$ in \cite{DLY} with
suitable modifications. 

We briefly explain how the decompostion of $V_{\sqrt{2}D_{l}}$
is obtained. We consider a sublattice $N$ of
$L=A_{1}^{\oplus l}$ isomorphic to $\sqrt{2}D_{l}$ in a standard
way. We show that there exists an
automorphism $\rho $ of the vertex operator algebra $V_{L}$ which maps 
$V_{N} $ onto $V_{L}^{+}$ (the fixed point subspace of the involution of 
$V_{L}$ induced from the involution $\beta \longmapsto -\beta $ of
$L$). Thus, the decomposition of $V_{N}$ as a $T$-module is
equivalent to the
decomposition of $V_{L}^{+}$ as a $\widetilde{T}$-module, where 
$\widetilde{T}=\rho (T)$. The decomposition of $V_{L}^{+}$ as 
a $\widetilde{T}$-module is relatively easy and  
 is completely determined in Section 4. For the purpose of completeness
we  also review the 
necessary results on the decomposition of a vertex operator algebra
associated with a positive definite rank one lattice as a module for 
$L(1,0)$ in Section 5, most of which are taken from \cite{DG} and \cite{KR}.

The authors would like to thank 
Masaaki Kitazume and Masahiko Miyamoto for helpful advice.

\section{Some automorphisms of $V_L$}

We shall use standard notation for the vertex operator algebra
$V_{L} = M(1) \otimes {\mathbb C}\{L\}$ associated with a positive
definite even lattice $L$ \cite{FLM}. Note that if the square
length of any element in $L$ is a
multiple of 4 or if $L$ is a rank one lattice then the central extension 
$\hat L$ of $L$ by the cyclic group of order two is split and 
${\mathbb C}\{L\}$ in the notation of \cite{FLM} is isomorphic to 
the group algebra ${\mathbb C}[L]$.

In \cite{DLY} three automorphisms $\theta_{1}$, $\theta_{2}$,
$\sigma$ of order two of the vertex operator algebra $V_{{\mathbb
Z}\alpha}$ associated with a rank one lattice ${\mathbb Z}\alpha$,
where $\langle \alpha, \alpha \rangle = 2$, are considered. They
are determined  by
\begin{align*}
\theta_{1} & : \alpha(-1) \longmapsto \alpha(-1), & e^{\alpha} & \longmapsto
-e^{\alpha}, & e^{-\alpha} & \longmapsto -e^{-\alpha}, \\
\theta_{2} & : \alpha(-1) \longmapsto -\alpha(-1), & e^{\alpha} &
\longmapsto e^{-\alpha}, & e^{-\alpha} & \longmapsto e^{\alpha}, \\
\sigma & : \alpha(-1) \longmapsto e^{\alpha} + e^{-\alpha}, & e^{\alpha} +
e^{-\alpha} & \longmapsto \alpha(-1), & e^{\alpha} - e^{-\alpha} &
\longmapsto -(e^{\alpha} - e^{-\alpha}).
\end{align*}

The automorphism $\theta_1$ maps $u \otimes e^{\beta}$ to 
$(-1)^{\langle \alpha, \beta \rangle/2} u
\otimes e^{\beta}$ for $u\in M(1)$ and $\beta \in {\mathbb Z}\al$ 
and $\theta_{2}$ is the automorphism induced from the isometry $\beta
\longmapsto -\beta$ of ${\mathbb Z}\alpha$. Note also that
\[
\sigma \theta_{1} \sigma = \theta_{2},\quad \sigma(\alpha(-1)^{2}) =
\alpha(-1)^{2}\text{\quad and \quad} \sigma(e^{\pm \alpha}) = (\alpha(-1)
\mp (e^{\alpha} - e^{-\alpha}))/2.
\]

Let $L$ be a lattice with basis $\{\alpha_{1},\alpha_{2},\ldots,
\alpha_{l}\} $ such that $\langle \alpha_{i},\alpha_{j}\rangle =
2\delta_{ij}$. Then, $L= {\mathbb Z}\alpha_1 \oplus \cdots \oplus {\mathbb Z}
\alpha_l$ and each ${\mathbb Z}\alpha_i$ is isomorphic to a root lattice of
type $A_1$. Since $\langle \alpha_i, \alpha_j \rangle = 0$ for $i \ne j$,
the vertex operator algebra $V_L$ is a tensor product 
$V_L = V_{{\mathbb Z}\alpha_1} \otimes \cdots \otimes 
V_{{\mathbb Z}\alpha_l}$ of $V_{{\mathbb Z}\alpha_i}$'s. 
Using the automorphisms $\theta_1$, $\theta_2$, and $\sigma$
of $V_{{\mathbb Z}\alpha_i}$ described above, we can define three 
automorphisms $\psi_1$, $\psi_2$, and $\tau$ of $V_L$ of order two by
\[
\psi_1 = \theta_1 \otimes \cdots \otimes \theta_1, \qquad \psi_2 = \theta_2
\otimes \cdots \otimes \theta_2, \qquad \tau = \sigma \otimes \cdots \otimes
\sigma.
\]
Then
\begin{equation}
\psi_{1}(u \otimes e^{\beta}) = (-1)^{\langle \alpha_1+\alpha_2+\cdots
+\alpha_l, \,\beta \rangle /2}u \otimes e^{\beta}
\end{equation}
for $u \in M(1)$ and $\beta \in L$, 
$\psi_2$ is the automorphism induced from the isometry $\beta \longmapsto
-\beta$ of $L$, and $\tau \psi_1 \tau = \psi_2$.

We shall calculate the images of some elements in $V_L$ under the
automorphism $\tau$, which will be used in the next section. If $i \ne j$,
then in $V_{{\mathbb Z}\alpha_i} \otimes V_{{\mathbb Z}\alpha_j}$ we have
\begin{align*}
\tau (\alpha_i(-1)\alpha_j(-1)) &= \sigma (\alpha_i(-1)) \sigma
(\alpha_j(-1)) \\
&= (e^{\alpha_i} + e^{-\alpha_i})(e^{\alpha_j} + e^{-\alpha_j}) \\
&= e^{\alpha_i+\alpha_j} + e^{\alpha_i-\alpha_j} + e^{-\alpha_i+\alpha_j} +
e^{-\alpha_i-\alpha_j}.
\end{align*}

Since $(\alpha_i\pm\alpha_j)(-1)^2 = \alpha_i(-1)^2 + \alpha_j(-1)^2 \pm
2\alpha_i(-1)\alpha_j(-1)$ in $V_{{\mathbb Z}\alpha_i} \otimes 
V_{{\mathbb Z}\alpha_j}$ and $\sigma (\alpha_i(-1)^2) = \alpha_i(-1)^2$, 
it follows that
\[
\tau((\alpha_i\pm\alpha_j)(-1)^2) = \alpha_i(-1)^2 + \alpha_j(-1)^2 \pm
2\tau(\alpha_i(-1)\alpha_j(-1)).
\]

We also have
\begin{align*}
\tau (e^{\alpha _{i}\pm \alpha _{j}})& =\sigma (e^{\alpha _{i}})\sigma
(e^{\pm \alpha _{j}}) \\
& =\frac{1}{4}(\alpha _{i}(-1)-(e^{\alpha _{i}}-e^{-\alpha _{i}}))(\alpha
_{j}(-1)\mp (e^{\alpha _{j}}-e^{-\alpha _{j}}))
\end{align*}
in $V_{{\mathbb Z}\alpha _{i}}\otimes V_{{\mathbb Z}\alpha _{j}}$. Likewise
\[
\tau (e^{-(\alpha _{i}\pm \alpha _{j})})=\frac{1}{4}(\alpha
_{i}(-1)+(e^{\alpha _{i}}-e^{-\alpha _{i}}))(\alpha _{j}(-1)\pm (e^{\alpha
_{j}}-e^{-\alpha _{j}})),
\]
and thus
\[
\tau (e^{\alpha _{i}\pm \alpha _{j}}+e^{-(\alpha _{i}\pm \alpha _{j})})=
\frac{1}{2}(\alpha _{i}(-1)\alpha _{j}(-1)\pm e^{\alpha _{i}+\alpha _{j}}\mp
e^{\alpha _{i}-\alpha _{j}}\mp e^{-\alpha _{i}+\alpha _{j}}\pm e^{-\alpha
_{i}-\alpha _{j}}).
\]

For $\beta = \pm (\alpha_i \pm \alpha_j)/\sqrt{2}$, set
\begin{equation}
w^{\pm}(\beta) = \frac{1}{2}\beta(-1)^2 \pm (e^{\sqrt{2}\beta} + 
e^{-\sqrt{2}\beta}).
\end{equation}
Then, $w^{\pm}(\beta) = w^{\pm}(-\beta)$. From the above
calculation, we have
\begin{equation}
\begin{split}
\ & \tau(w^-((\alpha_i \pm \alpha_j)/\sqrt{2})) = w^{\pm}((\alpha_i -
\alpha_j)/\sqrt{2}), \\
& \tau(w^+((\alpha_i \pm \alpha_j)/\sqrt{2})) = w^{\pm}((\alpha_i +
\alpha_j)/\sqrt{2}).
\end{split}
\end{equation}

\section{Conformal vectors}

Let us consider a sublattice
\[
N=\sum_{i,j=1}^{l}{\mathbb Z}(\alpha _{i}\pm \alpha _{j})
\]
of $L$, which is isomorphic to the root lattice of type $\sqrt{2}D_{l}$. We
choose the following elements as the simple roots of type $D_{l}$:
\begin{align*}
\beta _{1}& =(\alpha _{1}+\alpha _{2})/\sqrt{2},\qquad \beta _{2}=(-\alpha
_{2}+\alpha _{3})/\sqrt{2},\qquad \beta _{3}=
(-\alpha _{1}+\alpha _{2})/\sqrt{2}, \\
\beta _{i}& =(-\alpha _{i}+\alpha _{i+1})/\sqrt{2}\quad \mbox{ \qquad for } 
\quad 3\le i\le l-1.
\end{align*}
Then
\[
\Phi _{l}^{+}=\{(\alpha _{i}+\alpha _{j})/\sqrt{2},\,(-\alpha _{i}+\alpha
_{j})/\sqrt{2}\,|\,1\le i<j\le l\}
\]
is the set of positive roots. Using the notation $w^{\pm }(\beta )$ defined
in (2.2), we set
\begin{equation}
\begin{aligned}
s^{1}& =\frac{1}{4}w^{-}(\beta _{1}),\\
s^{2}& =\frac{1}{5}\left( w^{-}(\beta _{1})+w^{-}(\beta _{2})+w^{-}(\beta
_{1}+\beta _{2})\right),\\
s^{r}& =\frac{1}{2r}\sum_{1\le i<j\le r}\left( w^{-}\left( (\alpha
_{i}+\alpha _{j})/\sqrt{2}\right) +w^{-}\left( (-\alpha _{i}+\alpha _{j})/
\sqrt{2}\right) \right) ,\quad 3\le r\le l,\\
\omega & =\frac{1}{4(l-1)}\sum_{\beta \in \Phi _{l}^{+}}\beta (-1)^{2}.
\end{aligned}
\label{s} 
\end{equation}

It was shown by \cite{DLMN} that the elements
\begin{equation}
\omega ^{1}=s^{1},\qquad \omega ^{i}=s^{i}-s^{i-1},2\le i\le l,\qquad \omega
^{l+1}=\omega -s^{l}  \label{ea3}
\end{equation}
are mutually orthogonal conformal vectors. Their central charges $c(\omega
^{i})$ are as follows:
\[
c(\omega ^{1})=1/2,\quad c(\omega ^{2})=7/10,\quad c(\omega ^{3})=4/5,
\text{\quad and \quad }c(\omega ^{i})=1\text{ for }4\le i\le l+1.
\]
The subalgebra $\mbox{Vir}(\omega ^{i})$ of the vertex operator
algebra $V_{N}$ generated by $\omega ^{i}$ is isomorphic to the Virasoro 
vertex operator algebra $L(c(\omega ^{i}),0)$ which is the 
irreducible highest weight 
module for the Virasoro algebra with central charge $c(\omega^i)$ and 
highest weight $0.$ 
Since $\omega ^{1}$, 
$\omega^{2}$, \ldots , $\omega ^{l+1}$ are mutually orthogonal, 
the subalgebra $T$ 
of $V_{N}$ generated by these conformal vectors is a tensor product of 
$\mbox{Vir}(\omega ^{i})$'s, namely,
\begin{align*}
T& =\mbox{Vir}(\omega ^{1})\otimes \cdots \otimes \mbox{Vir}(\omega ^{l+1})
\\
\ & \cong L(c(\omega ^{1}),0)\otimes \cdots \otimes L(c(\omega ^{l+1}),0).
\end{align*}

As a $T$-module $V_N$ is completely reducible. Our purpose in this paper is
to determine all the irreducible direct summands of $V_{N}$ .

\begin{remark}
\label{choice} In \cite{DLMN} it is shown that corresponding to a chain 
$\cdots \subset \Phi ^{\prime \prime } \subset \Phi ^{\prime } 
\subset \Phi$ 
of root systems of type $A$, $D$, or $E$, one has a decomposition of the
Virasoro element of the vertex operator algebra $V_{\sqrt{2}R}$, where $R$
denotes a root lattice of type $\Phi $, into a sum of mutually orthogonal
conformal vectors (see also \cite[Lemma 5.1]{M}). In our notation $\Phi
_{l}=\Phi _{l}^{+}\cup -\Phi _{l}^{+}$ is a root system of type $D_{l}$ and 
$\omega$ is the Virasoro element of $V_{N}$. The conformal vectors $\omega
^{1}$, $\ldots $, $\omega ^{l+1}$ in (\ref{ea3}) correspond to the chain 
$\Phi _{1}\subset \Phi _{2}\subset \cdots \subset \Phi _{l}$, where $\Phi
_{1}=\{\pm \beta _{1}\}$, $\Phi _{2}=\{\pm \beta _{1},\pm \beta _{2},\pm
(\beta _{1}+\beta _{2})\}$, and
\[
\Phi _{r}=\{\pm (\alpha _{i}+\alpha _{j})/\sqrt{2},\,\pm (\alpha
_{i}-\alpha _{j})/\sqrt{2}\,|\,1\le i<j\le r\}, \quad 3\le r\le l.
\]
Note that $\Phi _{r}$ is a root system of type $A_{r}$ for $1\le r\le 3$ 
and of type $D_{r}$ for $4\le r\le l$.

There is another choice of conformal vectors. Namely, if one considers a
chain of root systems of type $A_{1}\subset A_{2}\subset \cdots \subset
A_{l-1}\subset D_{l}$, then the central charges of the corresponding
mutually orthogonal conformal vectors ${\omega ^{\prime }}^{1}$, $\ldots $, 
${\omega ^{\prime }}^{l+1}$ are given by
\[
c({\omega ^{\prime }}^{r})=1-\frac{6}{(r+2)(r+3)},\ 1\le r\le
l-1,\quad c({\omega ^{\prime }}^{l})=\frac{2(l-1)}{l+2},
\quad \text{and \quad }c({\omega ^{\prime }}^{l+1})=1.
\]
\end{remark}

The images of the conformal vectors $\omega^1$, $\ldots$, $\omega^{l+1}$
under the automorphism $\tau$ of $V_L$ can be calculated by (2.3). 
In fact, we have

\begin{lem}
\label{3.1} Let $s^{1},s^{2},\text{ and }s^{r},\ 3\leq r\leq l$ be defined 
as in (\ref{s}). Then,
\begin{align*}
\tau (s^{1})& =\frac{1}{4}w^{+}(\beta _{3}), \\
\tau (s^{2})& =\frac{1}{5}(w^{+}(\beta _{3})+w^{-}(\beta _{2})+
w^{+}(\beta_{2}+\beta _{3})), \\
\tau (s^{r})& =\frac{1}{2r}\sum_{1\le i<j\le r}(w^{+}((\alpha _{i}-\alpha
_{j})/\sqrt{2})+w^{-}((\alpha _{i}-\alpha _{j})/\sqrt{2})) \\
\ & =\frac{1}{4r}\sum_{1\le i<j\le r}(\alpha _{i}-\alpha _{j})(-1)^{2},\quad
3\le r\le l, \\
\tau (\omega )& =\omega .
\end{align*}
\end{lem}

\medskip 
Next we shall consider the images of these elements under the
automorphism $\varphi$ of $V_L$, where $\varphi$ is defined by 
\[
\varphi : u \otimes e^{\beta} \longmapsto (-1)^{\langle
\alpha_{2}+\alpha_{3}, \, \beta\rangle/2} u \otimes e^{\beta}
\]
for $u \in M(1)$ and $\beta \in L$. 
The automorphism $\varphi$ acts as $\theta_2$ on 
$V_{{\mathbb Z}\alpha_2}$ and $V_{{\mathbb Z}\alpha_3}$ and 
acts as the identity on 
$V_{{\mathbb Z}\alpha_i}$ for $i \ne 2$, $3$. Set $\rho = \varphi\tau$.

\begin{remark}
The automorphism $\rho $ here is slightly different from the automorphism 
$\rho $ of \cite{DLY}. In \cite{DLY} $\rho $ is defined to be a composite of 
$\theta _{2}\otimes 1\otimes 1$ and $\varphi \tau $. The symmetry $\alpha
_{1}\longleftrightarrow -\alpha _{1}$ and $\alpha _{i}\longleftrightarrow
\alpha _{i}$ for $2\le i\le l$ induces a ${\mathbb Z}_{2}$-symmetry of the
Dynkin diagram of type $D_{l}$, and the restriction of the automorphism 
$\theta _{2}\otimes 1\otimes \cdots \otimes 1$ of $V_{L}$ to $V_{N}$
corresponds to this symmetry. In this paper we do not consider it.
\end{remark}

The following lemma is easily obtained from Lemma \ref{3.1} and the
definition of $\varphi$.

\begin{lem}
We have
\begin{align*}
  \rho (s^{1})&=\frac{1}{4}w^{-}(\beta _{3}),& \rho (s^{2})&=\frac{1}{5}
  (w^{-}(\beta _{3})+w^{-}(\beta _{2})+w^{-}(\beta _{2}+\beta_{3})), \\
  \rho (s^{r})&=\tau (s^{r}),\quad 3\le r\le l,& \rho(\omega )&=\omega .
\end{align*}
\end{lem}

\medskip Let $\widetilde{\omega }^{i}=\rho (\omega ^{i})$ and set
\begin{equation}
\begin{aligned}
\gamma _{r}&=\alpha _{1}+\alpha _{2}+\cdots +\alpha _{r}-r\alpha
_{r+1},\quad 1\le r\le l-1, \\
\gamma _{l}&=\alpha _{1}+\alpha _{2}+\cdots +\alpha _{l}.
\end{aligned}
\end{equation}

\begin{lem}
\begin{description}
\item[$(1)$] The vectors $\widetilde{\omega }^{1}$, 
$\widetilde{\omega }^{2}$, and 
$\widetilde{\omega }^{3}$ are the mutually orthogonal conformal vectors of 
$V_{{\mathbb Z}(\alpha _{1}-\alpha _{2})+{\mathbb Z}(\alpha _{2}-\alpha
_{3})}\cong V_{\sqrt{2}A_{2}}$ defined in \cite{DLMN}.

\item[$(2)$] $\widetilde{\omega }^{r+1}={\displaystyle \frac{1}{4r(r+1)}\gamma
_{r}(-1)^{2}}$ \quad for $3\le r\le l-1$.

\item[$(3)$] $\widetilde{\omega }^{l+1}={\displaystyle \frac{1}{4l}\gamma
_{l}(-1)^{2}}$.
\end{description}
\end{lem}

\noindent {\bfseries Proof} \ (1) is shown in \cite[Lemma 3.6]{DLY}. For 
$3 \le r \le l-1$ we have
\begin{align*}
\rho(s^{r+1}) - \rho(s^r) &= \frac{1}{4(r+1)}\sum_{1 \le i < j \le
r+1}(\alpha_i - \alpha_j)(-1)^2 - \frac{1}{4r}\sum_{1 \le i < j \le
r}(\alpha_i - \alpha_j)(-1)^2 \\
\ &= \frac{1}{4r(r+1)}\Bigl( r\sum_{1 \le i \le r} (\alpha_i -
\alpha_{r+1})(-1)^2 - \sum_{1 \le i < j \le r}(\alpha_i -
\alpha_j)(-1)^2\Bigr) \\
\ &= \frac{1}{4r(r+1)}(\alpha_1+\alpha_2+\cdots+\alpha_r-r\alpha_{r+1})(-1)^2.
\end{align*}

We also have
\begin{align*}
\omega - \rho(s^l) &= \frac{1}{8(l-1)}\sum_{1 \le i < j \le l} ((\alpha_i +
\alpha_j)(-1)^2 + (\alpha_i - \alpha_j)(-1)^2) \\
\ & \qquad \qquad - \frac{1}{4l}\sum_{1 \le i < j \le l} (\alpha_i -
\alpha_j)(-1)^2 \\
\ &= \frac{1}{4(l-1)}\sum_{1 \le i < j \le l} (\alpha_i(-1)^2 +
\alpha_j(-1)^2) \\
\ & \qquad \qquad - \frac{1}{4l}\sum_{1 \le i < j \le l} (\alpha_i(-1)^2 +
\alpha_j(-1)^2 - 2\alpha_i(-1)\alpha_j(-1)) \\
\ &= \frac{1}{4l}(\alpha_1 + \alpha_2 + \cdots + \alpha_l)(-1)^2. \quad \qed
\end{align*}

Note that for $3 \le r \le l$, the element $\widetilde{\omega}^{r+1}$ is the
Virasoro element of the vertex operator algebra $V_{{\mathbb Z}\gamma_r}$
associated with a rank one lattice ${\mathbb Z}\gamma_r$.

Set $U^{\pm} = \{ v \in U \,|\,\psi_2(v) = \pm v\}$ for any 
$\psi_2$-invariant subspace $U$ of $V_L$.

\begin{lem}
$(1)$ $N=\{\beta \in L\,|\,\langle \alpha _{1}+\cdots +\alpha _{l},\beta
\rangle \equiv 0\pmod 4 \}$.

$(2)$ $V_{N}=\{v\in V_{L}\,|\,\psi _{1}(v)=v\}$.

$(3)$ $\rho (V_{N})=V_{L}^{+}$.
\end{lem}

\noindent {\bfseries Proof} \ (1) is clear from the definition of $N$. (1)
and (2.1) imply (2). Since $\tau \psi_1 = \psi_2\tau$ and the automorphisms 
$\psi_2$ and $\varphi$ commute, we have $\rho\psi_1 = \psi_2\rho$. Hence (3)
holds. \qed

\medskip The last assertion of the above lemma implies that the
decomposition of $V_N$ into a direct sum of irreducible $T$-modules is
equivalent to that of $V_L^+$ as a $\widetilde{T}$-module, where 
$\widetilde{T} = \rho (T)$ is of the form
\begin{align*}
\widetilde{T} &= \mbox{Vir}(\widetilde{\omega}^1) \otimes \cdots \otimes 
\mbox{Vir}(\widetilde{\omega}^{l+1}) \\
\ & \cong L(\frac{1}{2}, 0) \otimes L(\frac{7}{10}, 0) \otimes L(\frac{4}{5}
, 0) \otimes L(1, 0) \otimes \cdots \otimes L(1, 0).
\end{align*}

The decomposition of $V_L^+$ as a $\widetilde{T}$-module will be
accomplished in the next section.

\section{Decomposition of $V_L^+$ as a $\widetilde{T}$-module}

In this section, we shall study the decomposition of $V_{L}^{+}$ into a
direct sum of irreducible $\widetilde{T}$-modules. As in \cite{DLY} we set
\[
E={\mathbb Z}(\alpha _{1}-\alpha _{2})+{\mathbb Z}(\alpha _{2}-\alpha
_{3})\qquad \mbox{and}\qquad D=E+{\mathbb Z}\gamma _{3}+\cdots +{\mathbb Z}
\gamma _{l}.
\]
The elements $\alpha _{1}-\alpha _{2}$, $\alpha _{2}-\alpha _{3}$, $\gamma
_{3}$, $\ldots $, $\gamma _{l}$ form a basis of the lattice $D$. Since 
$\gamma _{r}=\sum_{i=1}^{r}i(\alpha _{i}-\alpha _{i+1})$, we can take
\[
\{\alpha _{1}-\alpha _{2},\,\alpha _{2}-\alpha _{3},\,3(\alpha
_{3}-\alpha _{4}),\,\ldots ,\,(l-2)(\alpha _{l-2}-\alpha
_{l-1}),\,\gamma _{l-1},\,\gamma _{l}\}
\]
as another basis. The lattices $E$, ${\mathbb Z}\gamma _{3}$, $\ldots $, 
${\mathbb Z}\gamma _{l}$ are mutually orthogonal, so the vertex operator
algebra $V_{D}$ associated with the lattice $D$ is a tensor product 
$V_{D}=V_{E}\otimes V_{{\mathbb Z}\gamma _{3}}\otimes \cdots \otimes 
V_{{\mathbb Z}\gamma _{l}}$.

Next, we want to describe the cosets of $D$ in $L$. Set
\[
\xi _{r}=\frac{1}{r(r+1)}\gamma _{r},\,1\le r\le l-1,\qquad \text{and}\qquad
\xi _{l}=\frac{1}{l}\gamma _{l}.
\]
Clearly $\langle \gamma _{r},\xi _{r}\rangle =2$ and we have
\begin{equation}
-\xi _{1}+\xi _{2}=\frac{1}{3}(-(\alpha _{1}-\alpha _{2})+(\alpha
_{2}-\alpha _{3})).
\end{equation}
Moreover,
\[
\alpha _{r}-\alpha _{r+1}=
\frac{1}{r}(-\gamma _{r-1}+\gamma_{r})=-(r-1)\xi _{r-1}+(r+1)\xi _{r},
\]
and so
\begin{equation}
-\xi _{1}+\xi _{2}+\cdots +\xi _{l}=\alpha _{2}.
\end{equation}

To simplify the notation, we set $\eta = -\xi_1 + \xi_2$.

\begin{lem}
$|D+{\mathbb Z}\alpha _{2}:D|$ is equal to the least common multiple of $3$,
$4$, $\ldots $, $l$.
\end{lem}

\noindent {\bfseries Proof} \ Since $\alpha_1 - \alpha_2$, 
$\alpha_2-\alpha_3 $, $\gamma_3$, $\ldots$, $\gamma_l$ are linearly
independent over ${\mathbb Q}$, the assertion follows from (4.1) and (4.2).
\qed

\medskip 
Lemma 4.1 is valid even for $l=3$ if we consider a lattice of type 
$A_{3}$ as of type $D_{3}$.

Note that $D + {\mathbb Z}\alpha_2 = L$ for $3 \le l \le 5$. Indeed, the
coset $D + \alpha_2$ contains $\alpha_1$, $\alpha_2$, and $\alpha_3$.
Moreover, $\alpha_4 \in D + 9\alpha_2$ if $l=4$, and $\alpha_4 \in D +
21\alpha_2$ and $\alpha_5 \in D + 36\alpha_2$ if $l=5$. Hence $n\alpha_2$, 
$0 \le n \le d-1$, where $d$ denotes the least common multiple of $3$, $4$, 
$\ldots$, $l$, form a complete system of representatives of the cosets of $D$
in $L$ in these three cases. However, $D + {\mathbb Z}\alpha_2 \ne L$ for 
$l\ge 6$. We shall use the following elements to describe all the cosets of $D$
in $L$. For $m_3,...,m_{l-2}, n\in \Z$ we let
\begin{equation}
\begin{split}
\lambda &= \lambda(m_3, m_4, \ldots , m_{l-2}, n) \\
&= m_3(\alpha_3-\alpha_4) + m_4(\alpha_4-\alpha_5) + \cdots +
m_{l-2}(\alpha_{l-2}-\alpha_{l-1}) + n\alpha_2 \\
&\equiv (m_3+n)\eta + \sum_{r=3}^{l-3} ((r+1)m_r - rm_{r+1} + n)\xi_r \\
& \qquad \quad +((l-1)m_{l-2}+n)\xi_{l-2} + n\xi_{l-1} + n\xi_l \quad \pmod D.
\end{split}
\end{equation}
The last congruence modulo $D$ comes from (4.1), (4.2), and the fact that 
$\alpha _{r}-\alpha _{r+1}=-(r-1)\xi _{r-1}+(r+1)\xi _{r}$.

\begin{lem}
\begin{description}
\item[$(1)$]  $\left\{ \lambda =\lambda (m_{3},\ldots ,m_{l-2},n)|\,0\le
m_{r}\le r-1,0\le n\le l(l-1)-1\right\} $ forms a complete system of
representatives of the cosets of $D$ in $L$.

\item[$(2)$]  Every element in the coset $D+\lambda $ can be uniquely
written in the form
\begin{align*}
(\nu +(m_{3}+n)\eta )& +\sum_{i=3}^{l-3}(\mu _{i}+((i+1)m_{i}-im_{i+1}+n)\xi
_{i}) \\
& +(\mu _{l-2}+((l-1)m_{l-2}+n)\xi _{l-2})+(\mu _{l-1}+n\xi _{l-1})+(\mu
_{l}+n\xi _{l})
\end{align*}
for $\nu \in E$ and $\mu _{i}\in {\mathbb Z}\gamma _{i}$.
\end{description}
\end{lem}

\noindent {\bfseries Proof} \ Let $M_r={\mathbb Z}(\alpha_1-\alpha_2) +
\cdots + {\mathbb Z}(\alpha_r-\alpha_{r+1}) +{\mathbb Z}\gamma_{r+1}+\cdots+
{\mathbb Z}\gamma_l$. Then $D = M_2 \subset M_3 \subset \cdots \subset 
M_{l-2} \subset L$. Note that
\begin{equation}
-l(\alpha_1-\alpha_2) + \sum_{i=2}^{l-2} l(l-i-1)(\alpha_i-\alpha_{i+1}) +
\gamma_{l-1} + (l-1)\gamma_l = l(l-1)\alpha_2.
\end{equation}
Since $\alpha _{i}\in M_{l-2}+\alpha _{2}$ for $1\le i\le l-1$
and $\alpha _{l}\in M_{l-2}+(l-1)^{2}\alpha _{2}$, 
$L=M_{l-2}+{\mathbb Z}\alpha _{2}$ and 
$\left\{ n\alpha _{2}|\,0\le n\le l(l-1)-1\right\} $ is a
complete system of representatives of the cosets of $M_{l-2}$ in $L$.
Furthermore, $\gamma _{r}=\sum_{i=1}^{r}i(\alpha _{i}-\alpha _{i+1})$
implies that $\left\{ m_{r}(\alpha _{r}-\alpha _{r+1})|0\le m_{r}\le
r-1\right\} $ is a complete system of representatives of the cosets of 
$M_{r-1}$ in $M_{r}$. Hence (1) holds. For $i=1$, $2$, the elements $\alpha
_{i}-\alpha _{i+1}$, $\gamma _{3}$, $\ldots $, $\gamma _{l}$ are mutually
orthogonal, and so each element of $D+\lambda $ is uniquely decomposed into
an orthogonal sum of the form stated in (2). \qed

\medskip By Lemma 4.2, we have $V_L = \oplus_{\lambda} V_{D+\lambda}$ where
the sum runs through the coset representatives of $D$ in $L,$  
and
\begin{equation}
\begin{split}
V_{D+\lambda} &= V_{E+(m_3+n)\eta} \otimes \Bigl( \otimes_{r=3}^{l-3} 
V_{{\mathbb Z}\gamma_r +((r+1)m_r- rm_{r+1}+n)\xi_r} \Bigr) \\
& \qquad \otimes V_{{\mathbb Z}\gamma_{l-2}+((l-1)m_{l-2}+n)\xi_{l-2}}
\otimes V_{{\mathbb Z}\gamma_{l-1}+n\xi_{l-1}} \otimes 
V_{{\mathbb Z}\gamma_l + n\xi_l}
\end{split}
\end{equation}
for $\lambda = \lambda(m_3, \ldots , m_{l-2}, n)$. Since $\langle E,
\eta\rangle = 2{\mathbb Z}$, the tensor factor of the form $V_{E+k\eta}$ is
an irreducible module for $V_E$ \cite{D}. Likewise, 
$V_{{\mathbb Z}\gamma_i + k\xi_i}$ is an irreducible module 
for $V_{{\mathbb Z}\gamma_i}$.
Thus (4.5) holds as a module for $V_D = V_E \otimes V_{{\mathbb Z}\gamma_3}
\otimes \cdots \otimes V_{{\mathbb Z}\gamma_l}$.

\begin{lem}
For $\lambda =\lambda (m_{3},\ldots ,m_{l-2},n)$, $0\le m_{r}\le r-1$, 
$0\le n\le l(l-1)-1$, we have $D+\lambda =D-\lambda $ if and only if $m_{3}$,
$\ldots $, $m_{l-2}$, and $n$ satisfy one of the following conditions.

\begin{description}
\item[$(1)$]  $n=0$, and $m_{r}=0$ if $r$ is odd and $m_{r}=0$ or $r/2$ if $r$
is even.

\item[$(2)$]  $n=l(l-1)/2$ and $2m_{r}+l(l-1)\equiv 0\pmod
r$. Such an $m_{r}$ is unique if $r$ is odd and there are exactly two such 
$m_{r}$ if $r$ is even.
\end{description}
\end{lem}

\noindent {\bfseries Proof} \ Suppose $D+\lambda = D-\lambda$, or
equivalently $2\lambda \in D$. Then $2n\alpha_2 \in M_{l-2}$ since $\lambda$
and $n\alpha_2$ are congruent modulo $M_{l-2}$. Now the quotient group 
$L/M_{l-2}$ is a cyclic group of order $l(l-1)$ generated by the coset 
$M_{l-2}+\alpha_2$. Thus $n=0$ or $l(l-1)/2$.

First, we assume that $n=0$. In this case, since $\alpha_1-\alpha_2$, 
$\alpha_2-\alpha_3$, $3(\alpha_3-\alpha_4)$, $\ldots$, $(l-2)(\alpha_{l-2}-
\alpha_{l-1})$, $\gamma_{l-1}$, $\gamma_l$ form a basis of the lattice $D$,
we have $2\lambda \in D$ if and only if $2m_r$ is divisible by $r$ for 
$3\le r \le l-2$.

Next, assume that $n = l(l-1)/2$. In this case
\begin{align*}
2\lambda &= 2m_3(\alpha_3-\alpha_4) + \cdots +
2m_{l-2}(\alpha_{l-2}-\alpha_{l-1}) + l(l-1)\alpha_2 \\
&= -l(\alpha_1-\alpha_2) + l(l-3)(\alpha_2-\alpha_3) \\
& \qquad + \sum_{r=3}^{l-2}(2m_r+l(l-r-1))(\alpha_r-\alpha_{r+1}) +
\gamma_{l-1} + (l-1)\gamma_l
\end{align*}
by (4.4). Hence $2\lambda \in D$ if and only if $2m_r + l(l-r-1)$ is
divisible by $r$ for $3 \le r \le l-2$. \qed

\medskip The automorphism $\psi_2$ fixes the conformal vectors 
$\widetilde{\omega}^1$, $\ldots$, $\widetilde{\omega}^{l+1}$, 
and so $\widetilde{T}
\subset V_D^+$. In particular, $\psi_2$ is a $\widetilde{T}$-module
isomorphism. We have $\psi_2(V_{D+\lambda}) = V_{D-\lambda}$, and thus 
$V_{D-\lambda}$ is isomorphic to $V_{D+\lambda}$ as a $\widetilde{T}$-module.
If $D+\lambda \ne D-\lambda$, the fixed point subspace $(V_{D+\lambda}
\oplus V_{D-\lambda})^+$ in $V_{D+\lambda} \oplus V_{D-\lambda}$ is equal to
$\{ v+ \psi_2(v)\,|\, v \in V_{D+\lambda}\}$ and it is isomorphic to 
$V_{D+\lambda}$.

If $D+\lambda = D-\lambda$ for $\lambda = \lambda(m_3, \ldots , m_{l-2}, n)$, 
then $m_3$, $\ldots$, $m_{l-2}$, and $n$ satisfy the conditions in Lemma
4.3. In this case (4.5) is in the following form:
\[
V_{D+\lambda} = V_E \otimes V_{{\mathbb Z}\gamma_3+b_3\gamma_3} \otimes
\cdots \otimes V_{{\mathbb Z}\gamma_l+b_l\gamma_l},
\]
with $b_i \in \{0, \frac{1}{2}\}$ for $3 \le i \le l-1$ and $b_l=0$ or $b_l
\in \{0, \frac{1}{2}\}$ depending on whether $l$ is odd or even.

For $\varepsilon = (\varepsilon_0, \varepsilon_3, \ldots,
\varepsilon_l)$ with $\varepsilon_i = +$ or $-$, set
\begin{equation}
V_{D+\lambda}^{\varepsilon} = V_E^{\varepsilon_0} \otimes 
V_{{\mathbb Z}\gamma_3+b_3\gamma_3}^{\varepsilon_3} \otimes \cdots 
\otimes V_{{\mathbb Z}\gamma_l+b_l\gamma_l}^{\varepsilon_l}.
\end{equation}
Then
\begin{equation}
V_{D+\lambda}^+ = \oplus_{\varepsilon} V_{D+\lambda}^{\varepsilon},
\end{equation}
where $\varepsilon$ runs over all $\varepsilon = (\varepsilon_0,
\varepsilon_3, \ldots, \varepsilon_l)$ such that even number of 
$\varepsilon_i$'s are $-$.

The decomposition of $V_{E}^{\pm }$ and $V_{E+\eta }$ as 
$\widetilde{T}^{\prime }$-modules, where $\widetilde{T}^{\prime }=
\mbox{Vir}(\widetilde{\omega }^{1})\otimes 
\mbox{Vir}(\widetilde{\omega }^{2})\otimes 
\mbox{Vir}(\widetilde{\omega }^{3})$, can be found in 
\cite[Lemma 4.1]{KMY} and
\cite[(3.2), (3.3)]{DLY}. In fact, $V_{E}^{+}$ is a direct sum of four
irreducible $\widetilde{T}^{\prime }$-modules, which are isomorphic to
\begin{equation}
\begin{array}{ll}
\medskip L(\frac{1}{2},\,0)\otimes L(\frac{7}{10},\,0)\otimes 
L(\frac{4}{5},\,0),\qquad & L(\frac{1}{2},\,0)\otimes 
L(\frac{7}{10},\,\frac{3}{5})\otimes L(\frac{4}{5},\,\frac{7}{5}), \\
\medskip L(\frac{1}{2},\,\frac{1}{2})\otimes L(\frac{7}{10},\,
\frac{1}{10})\otimes L(\frac{4}{5},\,\frac{7}{5}), 
& L(\frac{1}{2},\,\frac{1}{2})\otimes
L(\frac{7}{10},\,\frac{3}{2})\otimes L(\frac{4}{5},\,0),
\end{array}
\end{equation}
and $V_{E}^{-}$ is a direct sum of four irreducible 
$\widetilde{T}^{\prime }$-modules, which are isomorphic to
\begin{equation}
\begin{array}{ll}
\medskip L(\frac{1}{2},\,0)\otimes L(\frac{7}{10},\,\frac{3}{5})\otimes 
L(\frac{4}{5},\,\frac{2}{5}),\qquad & L(\frac{1}{2},\,\frac{1}{2})\otimes 
L(\frac{7}{10},\,\frac{1}{10})\otimes L(\frac{4}{5},\,\frac{2}{5}), \\
\medskip L(\frac{1}{2},\,0)\otimes L(\frac{7}{10},\,0)\otimes 
L(\frac{4}{5},\,3), & L(\frac{1}{2},\,\frac{1}{2})\otimes 
L(\frac{7}{10},\,\frac{3}{2})\otimes L(\frac{4}{5},\,3).
\end{array}
\end{equation}
Similarly, $V_{E+\eta }$ is a direct sum of four irreducible 
$\widetilde{T}^{\prime }$-modules, which are isomorphic to
\begin{equation}
\begin{array}{ll}
\medskip L(\frac{1}{2},\,0)\otimes L(\frac{7}{10},\,0)\otimes 
L(\frac{4}{5},\,\frac{2}{3}),\qquad & L(\frac{1}{2},\,0)\otimes 
L(\frac{7}{10},\,\frac{3}{5})\otimes L(\frac{4}{5},\,\frac{1}{15}), \\
\medskip L(\frac{1}{2},\,\frac{1}{2})\otimes 
L(\frac{7}{10},\,\frac{1}{10})\otimes L(\frac{4}{5},\,\frac{1}{15}), 
& L(\frac{1}{2},\,\frac{1}{2})\otimes L(\frac{7}{10},\,\frac{3}{2})\otimes 
L(\frac{4}{5},\,\frac{2}{3}).
\end{array}
\end{equation}

For $3 \le i \le l$, the decomposition of $V_{{\mathbb Z}\gamma_i}^{\pm}$
and $V_{{\mathbb Z}\gamma_i + \gamma_i/2}$ as 
$\mbox{Vir}(\widetilde{\omega}^{i+1})$-modules is given 
in Appendix (see also \cite{DG}).

We divide a complete system of representatives of the cosets of $D$ in $L$
into three subsets $\Lambda_1$, $\Lambda_2$, and $-\Lambda_2$ so that $D +
\lambda = D - \lambda$ if and only if $\lambda \in \Lambda_1$. Then
\[
V_L = \Bigl( \oplus_{\lambda \in \Lambda_1} V_{D+\lambda} \Bigr) \oplus
\Bigl( \oplus_{\lambda \in \Lambda_2} V_{D+\lambda} \Bigr) \oplus \Bigl(
\oplus_{\lambda \in \Lambda_2} V_{D-\lambda} \Bigr).
\]

By the above argument, we conclude that

\begin{thm}
As $\widetilde{T}$-modules,
\[
V_{L}^{+}\cong \Bigl(\oplus _{\lambda \in \Lambda _{1}}V_{D+\lambda }^{+}
\Bigr)\oplus \Bigl(\oplus _{\lambda \in \Lambda _{2}}V_{D+\lambda }\Bigr).
\]
Furthermore, the decomposition of $V_{D+\lambda }^{+}$, 
$\lambda \in \Lambda _{1}$,
and $V_{D+\lambda }$, $\lambda \in \Lambda _{2}$, into a direct sum of
irreducible $\widetilde{T}$-modules is given by (4.5) through (4.10) and
(5.1) through (5.7).
\end{thm}

\medskip By the same method, we also have the decomposition of 
$V_{L}^{-}\cong (\oplus _{\lambda \in \Lambda _{1}}V_{D+\lambda }^{-})\oplus
(\oplus _{\lambda \in \Lambda _{2}}V_{D+\lambda })$ into a direct sum of
irreducible $\widetilde{T}$-modules.

\medskip 
We illustrate the first two cases, namely $l=3$ and $4$, and leave
the details for the other cases to the interested reader. In case of $l=3$, 
$N=\sqrt{2}A_{3}$ and the decomposition of $V_{L}^{+}$ as a 
$\widetilde{T}$-module given by Theorem 4.4 is identical with the 
decomposition of $V_{\sqrt{2}A_{3}}$ as a $T$-module obtained in 
\cite[Theorem 3.7]{DLY}.

If $l=4$, then $D=E+{\mathbb Z}\gamma _{3}+{\mathbb Z}\gamma
_{4}$, where $\gamma _{3}=\alpha _{1}+\alpha _{2}+\alpha
_{3}-3\alpha _{4}$ and $\gamma
_{4}=\alpha _{1}+\alpha _{2}+\alpha _{3}+\alpha _{4}$. By Lemma 4.2, 
$\{\lambda =\lambda (n)=n\alpha _{2}|\,0\le n\le 11\}$ is a
complete system of coset representatives of $D$ in $L$ with
$\Lambda _{1}=\{0,6\alpha _{2}\}$ and $\Lambda _{2}=
\{j\alpha_{2}\,|\,1\le j\le 5\}$. Thus
\[
V_{L}^{+}\cong V_{D}^{+}\oplus V_{D+6\alpha _{2}}^{+}\oplus (\oplus
_{j=1}^{5}V_{D+j\alpha _{2}}).
\]

Now, $\alpha _{2}=\eta +\xi _{3}+\xi _{4}$ and $V_{D+j\alpha
_{2}}=V_{E+j\eta }\otimes V_{{\mathbb Z}\gamma _{3}+j\xi _{3}}\otimes 
V_{{\mathbb Z}\gamma _{4}+j\xi _{4}}$. Furthermore,
\begin{eqnarray*}
V_{D}^{+} &=&(V_{E}^{+}\otimes V_{{\mathbb Z}\gamma _{3}}^{+}\otimes 
V_{{\mathbb Z}\gamma _{4}}^{+})\oplus (V_{E}^{+}\otimes V_{{\mathbb Z}\gamma
_{3}}^{-}\otimes V_{{\mathbb Z}\gamma _{4}}^{-}) \\
&&\oplus (V_{E}^{-}\otimes V_{{\mathbb Z}\gamma _{3}}^{+}\otimes 
V_{{\mathbb Z}\gamma _{4}}^{-})\oplus (V_{E}^{-}\otimes V_{{\mathbb Z}\gamma
_{3}}^{-}\otimes V_{{\mathbb Z}\gamma _{4}}^{+}),
\end{eqnarray*}
\begin{eqnarray*}
V_{D+6\alpha _{2}}^{+} &=&(V_{E}^{+}\otimes V_{{\mathbb Z}\gamma _{3}+\gamma
_{3}/2}^{+}\otimes V_{{\mathbb Z}\gamma _{4}+\gamma _{4}/2}^{+})\oplus
(V_{E}^{+}\otimes V_{{\mathbb Z}\gamma _{3}+\gamma _{3}/2}^{-}\otimes 
V_{{\mathbb Z}\gamma _{4}+\gamma _{4}/2}^{-}) \\
&&\oplus (V_{E}^{-}\otimes V_{{\mathbb Z}\gamma _{3}+\gamma
_{3}/2}^{+}\otimes V_{{\mathbb Z}\gamma _{4}+\gamma _{4}/2}^{-})\oplus
(V_{E}^{-}\otimes V_{{\mathbb Z}\gamma _{3}+\gamma _{3}/2}^{-}\otimes 
V_{{\mathbb Z}\gamma _{4}+\gamma _{4}/2}^{+}).
\end{eqnarray*}

Hence the decomposition of $V_{L}^{+}$ into a direct sum of irreducible 
$\widetilde{T}$-modules follows from (4.8) through (4.10) and (5.1) through
(5.7).

\section{Appendix}

We review the decomposition of a vertex operator algebra associated with a
rank one lattice and its modules into a direct sum of irreducible 
$L(1,0)$-modules. Most materials discussed here can be found in
\cite{DG}, \cite{KR}.

Let ${\mathbb Z}\gamma $ be a rank one positive definite even lattice with 
$\langle \gamma ,\gamma \rangle =2n$ and $n$ a positive integer.
Set $\xi =\gamma /2n$. Then $V_{{\mathbb Z}\gamma +a\xi }$, 
$0\le a\le 2n-1$, are the inequivalent irreducible
modules for the vertex
operator algebra $V_{{\mathbb Z}\gamma }$ \cite{D}. We remark that 
${\mathbb Z}\gamma $ is denoted by $L_{2n}$ in \cite{DG}. The element $\omega
=\gamma (-1)^{2}/4n$ is the Virasoro element of $V_{{\mathbb Z}\gamma }$.
The subalgebra $\mbox{Vir}(\omega )$ generated by $\omega $ is isomorphic to
the Virasoro vertex operator algebra $L(1,0)$ of central charge $1$. Note
that
\[
V_{{\mathbb Z}\gamma +a\xi }=\oplus _{m\in {\mathbb Z}}M(1)\otimes
e^{(m+a/2n)\gamma },
\]
and each direct summand $M(1)\otimes e^{(m+a/2n)\gamma }$ is completely
reducible as a $\mbox{Vir}(\omega )$-module (\cite[Proposition 3.1]{KR}).

If $W$ is a $\mbox{Vir}(\omega )$-module, the $q$-dimension of $W$ is
defined by $\dim _{q}W=\sum_{n}(\dim W_{n})q^{n}$. Let $L(1,h)$ be the
irreducible highest weight module for $L(1,0)$ with highest weight $h\in 
{\mathbb C}$. Its $q$-dimension is $\dim
_{q}L(1,h)=(q^{m^{2}/4}-q^{(m+2)^{2}/4})/\phi (q)$ if $h=m^{2}/4$ for some
integer $m$ and otherwise $\dim _{q}L(1,h)=q^{h}/\phi (q)$, where $\phi
(q)=\prod_{j=1}^{\infty }(1-q^{j})$ (see \cite{DG}). On the other hand, the
weight of $e^{(m+a/2n)\gamma }$ is
\[
\frac{1}{2}\langle (m+\frac{a}{2n})\gamma ,(m+\frac{a}{2n})\gamma \rangle = 
\frac{1}{4n}(2mn+a)^{2}.
\]
Then $\dim _{q}M(1)\otimes e^{(m+a/2n)\gamma }=q^{(2mn+a)^{2}/4n}/\phi (q)$.
Comparing this with $\dim _{q}L(1,h)$, we know how $M(1)\otimes
e^{(m+a/2n)\gamma }$ decomposes into a direct sum of irreducible modules 
$L(1,h)$'s.

As usual we consider the automorphism $\theta$ of order two induced from the
automorphism $\beta \longmapsto -\beta$ of the underlying lattice. Since 
$\theta$ fixes $\omega$, it is in fact an automorphism for $\mbox{Vir}
(\omega) $-modules. In particular, $\theta (M(1) \otimes e^{k\gamma}) = M(1)
\otimes e^{-k\gamma}$ is isomorphic to $M(1) \otimes e^{k\gamma}$ as 
$\mbox{Vir}(\omega)$-modules. If $U$ is a $\theta$-invariant
subspace, the eigenspace with eigenvalue $\pm 1$ is denoted by
$U^{\pm}$. Then $U = U^+ \oplus U^-$.

Note that $(2mn+a)^{2}/4n=j^{2}/4$ for some integer $j$ if and only if 
$n=k^{2}$ for some nonnegative integer $k$ and $a=0$, $n$. We divide the
calculation into three cases.

\medskip Case 1: $n=k^2$ for some $0 < k \in {\mathbb Z}$ and $a=0$. In this
case the weight of $e^{(m+a/2n)\gamma} = e^{m\gamma}$ is $m^2k^2$ and
\[
M(1) \otimes e^{m\gamma} = \oplus_{p=0}^{\infty} L(1, (mk+p)^2).
\]
As a special case, we have $M(1) = \oplus_{p=0}^{\infty} L(1,
p^2)$ when $m=0$. Moreover, $\theta (M(1)\otimes e^{m\gamma}) =
M(1)\otimes e^{-m\gamma}$ is isomorphic to $M(1)\otimes
e^{m\gamma}$ as a $\mbox{Vir}(\omega)$-module and
\[
V_{{\mathbb Z}\gamma} = \Bigl( \oplus_{m=1}^{\infty} M(1)\otimes
e^{-m\gamma}\Bigr) \oplus M(1) \oplus \Bigl( \oplus_{m=1}^{\infty}
M(1)\otimes e^{m\gamma}\Bigr),
\]
and we have (see \cite{DG})
\begin{align}
V_{{\mathbb Z}\gamma}^+ &\cong \Bigl( \oplus_{p=0}^{\infty} L(1, 4p^2)\Bigr)
\oplus \Bigl( \oplus_{m=1}^{\infty} \oplus_{p=0}^{k-1} mL(1, (mk+p)^2)\Bigr),
\\
V_{{\mathbb Z}\gamma}^- &\cong \Bigl( \oplus_{p=0}^{\infty} L(1,
(2p+1)^2)\Bigr) \oplus \Bigl( \oplus_{m=1}^{\infty} \oplus_{p=0}^{k-1} mL(1,
(mk+p)^2)\Bigr).
\end{align}

Case 2: $n=k^{2}$ for some $0<k\in {\mathbb Z}$ and $a=n$. The weight of 
$e^{(m+a/2n)\gamma }=e^{(m+1/2)\gamma }$ is $(2m+1)^{2}k^{2}/4$ and
\[
q^{(2m+1)^{2}k^{2}/4}=\sum_{p=0}^{\infty
}(q^{((2m+1)k+2p)^{2}/4}-q^{((2m+1)k+2p+2)^{2}/4})
\]
in the algebra ${\mathbb C}[[q^{1/4}]]$ of formal power series. Comparing
the $q$-dimension, we have
\[
M(1)\otimes e^{(m+1/2)\gamma }\cong \oplus _{p=0}^{\infty }
L(1,\frac{1}{4}((2m+1)k+2p)^{2}).
\]
Since $\theta (M(1)\otimes e^{(m+1/2)\gamma })=M(1)\otimes
e^{(-m-1+1/2)\gamma }$ and
\[
V_{{\mathbb Z}\gamma +\gamma /2}=\Bigl(\oplus _{m=0}^{\infty }M(1)\otimes
e^{(m+1/2)\gamma }\Bigr)\oplus \Bigl(\oplus _{m=0}^{\infty }M(1)\otimes
e^{(-m-1+1/2)\gamma }\Bigr),
\]
it follows that
\begin{equation}
\begin{split}
V_{{\mathbb Z}\gamma +\gamma /2}^{+}& \cong V_{{\mathbb Z}\gamma +\gamma
/2}^{-}\cong \oplus _{m=0}^{\infty }M(1)\otimes e^{(m+1/2)\gamma } \\
& \cong \oplus _{m=0}^{\infty }\oplus _{p=0}^{\infty }
L(1,\frac{1}{4}((2m+1)k+2p)^{2}).
\end{split}
\end{equation}

Case 3: $(2mn+a)^{2}/n\ne j^{2}$ for any integer $j$. In this case
\[
M(1)\otimes e^{(m+a/2n)\gamma }\cong L(1,\frac{1}{4n}(2mn+a)^{2})
\]
is an irreducible $\mbox{Vir}(\omega )$-module. Now $V_{{\mathbb Z}\gamma
+a\xi }$ is invariant under $\theta $ if and only if $a=0$ or $n$. If $a=0$,
then $M(1)\otimes e^{m\gamma }\cong L(1,m^{2}n)$ for $m\ne 0$ and
\begin{align}
V_{{\mathbb Z}\gamma }^{+}& \cong \Bigl(\oplus _{p=0}^{\infty }L(1,4p^{2})
\Bigr)\oplus \Bigl(\oplus _{m=1}^{\infty }L(1,m^{2}n)\Bigr), \\
V_{{\mathbb Z}\gamma }^{-}& \cong \Bigl(\oplus _{p=0}^{\infty}
L(1,(2p+1)^{2})\Bigr)\oplus \Bigl(\oplus _{m=1}^{\infty }L(1,m^{2}n)\Bigr)
\end{align}
(see \cite[Lemma 2.10]{DG}). If $a=n$, then we have
\begin{equation}
V_{{\mathbb Z}\gamma +\gamma /2}^{+}\cong V_{{\mathbb Z}\gamma +\gamma
/2}^{-}\cong \oplus _{m=0}^{\infty }L(1,\frac{1}{16}n(4m+1)^{2})
\end{equation}
as $\mbox{Vir}(\omega )$-modules. Finally, if $a\ne 0$, $n$, then
\begin{equation}
V_{{\mathbb Z}\gamma +a\xi }\cong \oplus _{m\in {\mathbb Z}}
L(1,\frac{1}{4n}(2mn+a)^{2}).
\end{equation}

\end{document}